\documentclass[reqno]{amsart}
\usepackage{amsfonts,amsmath,amssymb}

\newcommand{\beq}{\begin{equation}}
\newcommand{\eeq}{\end{equation}}
\newcommand{\bs}{\begin{split}}
\newcommand{\esplit}{\end{split}}
\newcommand{\begincal}{\begin{eqnarray*}}
\newcommand{\fincal}{\end{eqnarray*}}

\newtheorem{thm}{Theorem}[section]

\newtheorem{claim}{Claim}[section]

\newcommand{\eps}{\varepsilon}
\newcommand{\rn}{{\mathbb R}^n}

\title[Prescribing the scalar curvature in the null case]
{Prescribing the scalar curvature \\  in the null case}

\author{Olivier Druet}
\address{Olivier Druet \\
UMPA - ENS Lyon\\
46, all\'ee d'Italie\\
F-69364 Lyon Cedex 7\\
France}
\email{odruet@umpa.ens-lyon.fr}

\date{July, 2005}

\begin{document}

\begin{abstract} 
In this paper, we give a proof of the Kazdan-Warner conjecture concerning the prescribed scalar curvature problem in the null case.
\end{abstract}

\maketitle

\section{Introduction}

Let $(M,g)$ be a smooth compact Riemannian manifold of dimension $n\ge 3$. 
The Yamabe invariant of $\left(M,g\right)$ is a conformal invariant defined by 
\begin{equation}
\label{eqdefYamabeinvariant}
\mu_g = \frac{n-2}{4(n-1)}\inf_{\tilde{g}\in \left[g\right]} Vol_{\tilde{g}}\left(M\right)^{-\frac{n-2}{n}} \int_M S_{\tilde{g}}\, dv_{\tilde{g}}
\end{equation}
where 
$$\left[g\right]=\left\{e^u\, g,\, u\in C^\infty\left(M\right)\right\}$$
is the conformal class of $g$, $Vol_{\tilde{g}}\left(M\right)$ denotes the volume of $M$ with respect to the Riemannian measure associated to $\tilde{g}$ and $S_{\tilde{g}}$ denotes the scalar curvature 
of $\left(M,\tilde{g}\right)$. If $\tilde{g}=u^{\frac{4}{n-2}}g$ for some smooth positive 
function $u$, then the scalar curvatures of $g$ and $\tilde{g}$ are related by the following equation~:
\begin{equation}
\label{eqcourburescalaire}
\Delta_g u +\frac{n-2}{4(n-1)}S_g u= \frac{n-2}{4(n-1)}S_{\tilde{g}}u^{2^\star-1}
\end{equation}
where $2^\star=\frac{2n}{n-2}$ and $\Delta_g u = -div_g\left(\nabla u\right)$. This permits to rewrite 
(\ref{eqdefYamabeinvariant}) as 
\begin{equation}
\label{eqYamabeinvariant}
\mu_g = \inf_{u\in C^\infty\left(M\right),\, u\not\equiv 0} \frac{\int_M \left(\left\vert \nabla u\right\vert_g^2 
+\frac{n-2}{4(n-1)}S_g u^2\right)\, dv_g}{\left(\int_M \vert u\vert^{2^\star}\, dv_g\right)^{\frac{2}{2^\star}}}
\hskip.1cm.
\end{equation}
In this paper, we deal with the so-called prescribed scalar curvature problem in the null case. 
We assume that $\mu_g=0$ and we consider the following question~: given $f\in C^\infty\left(M\right)$, does there exist $\tilde{g}\in \left[g\right]$ such that $S_{\tilde{g}}=f$~? 
Thanks to equation (\ref{eqcourburescalaire}), and up to some obvious rescaling, the question is equivalent to~: given $f\in C^\infty\left(M\right)$, does there exist some smooth positive solution 
$u$ of 
$$\Delta_g u + \frac{n-2}{4(n-1)}S_g u = f u^{2^\star-1}\hskip.1cm?$$
The case $f\equiv 0$ corresponds to the well-known Yamabe problem in the null case, see Yamabe \cite{Yamabe}, and has been solved in this situation by Trudinger \cite{Trudinger}. Thus we 
may assume without loss of generality that $S_g\equiv 0$ and that $Vol_g\left(M\right)=1$ since 
there always exists such a metric in the conformal class of $g$. Then the above equation becomes 
\begin{equation}\label{eqprescribedscalarcurvature}
\Delta_g u = fu^{2^\star-1}\hskip.1cm.
\end{equation}
Obvious necessary conditions on $f$ so that there exists a smooth positive solution $u$ of the above equation are that 
\begin{equation}
\label{eqNC1}
\int_M f\, dv_g <0
\end{equation}
and that 
\begin{equation}
\label{eqNC2}
\max_M f>0
\end{equation}
if $f\not \equiv 0$, what we will assume in the sequel since the case $f\equiv 0$ is already settled. 
In order to prove that (\ref{eqNC1}) is a necessary condition, one has just to multiply equation (\ref{eqprescribedscalarcurvature}) by $u^{1-2^\star}$ and to integrate over $M$ to obtain that
$$\int_M f\, dv_g = \int_M u^{1-2^\star}\Delta_g u \, dv_g = -\frac{n+2}{n-2}\int_M u^{-2^\star}\left\vert \nabla u\right\vert_g^2\, dv_g <0$$
since $u\not \equiv Cst$ if $f\not\equiv 0$. As for the second necessary condition (\ref{eqNC2}), one has just to integrate equation (\ref{eqprescribedscalarcurvature}) over $M$ to obtain that 
$$\int_M f u^{2^\star-1}\, dv_g=0$$
which enforces $f$ to change sign if it is not identically zero. 

\medskip The natural question then is~: are the conditions (\ref{eqNC1}) and (\ref{eqNC2}) sufficient 
for $f$ to be the scalar curvature of a metric conformal to $g$~? In the sequel of their $2$-dimensional work \cite{KazdanWarner2}, Kazdan and Warner conjectured that the answer to this question was positive (see \cite{Kazdan, KazdanWarner1}). Jung, in \cite{Jung}, 
used the technique of sub- and upper-solution to solve equation (\ref{eqprescribedscalarcurvature}) under the assumptions (\ref{eqNC1}) and (\ref{eqNC2}). However, there is some problem in this proof - see p. 743 of \cite{Jung} - since the limit, for $s>0$ real, of 
$$\left(\frac{1-x^x}{s^{x^2}-x^x}\right)^{-\frac{\ln x}{x}}$$
as $x\to 0^+$ is not $1$ but $s^{-1}$. We provide in this paper a proof of the conjecture of Kazdan and Warner, based on fine blow-up analysis. In other words, we prove the following result~: 

\begin{thm}
\label{mainthm}
Let $(M,g)$ be a smooth compact Riemannian manifold of dimension $n\ge 3$ with null scalar curvature. 
A smooth function $f$ on $M$ is the scalar curvature of a metric conformal to $g$ if and only if 
$f\equiv 0$ or (\ref{eqNC1}) and (\ref{eqNC2}) hold.
\end{thm}

\medskip Note that this result had already been proved in dimensions $3$ and $4$ by minimisation techniques and test functions computations by Escobar and Schoen \cite{EscobarSchoen}. In this 
work, they also proved the existence of a solution of the equation (\ref{eqprescribedscalarcurvature}) 
under some flatness assumptions on $f$ around one of its maximum point in the case $(M,g)$ is locally 
conformally flat. 

\medskip The solution we obtain during the proof of the theorem 
is a solution of minimal energy, a so-called minimizing solution. 

\medskip Let us at last remark that, given $(M,g)$ a smooth compact Riemannian manifold of dimension $n\ge 3$ with null Yamabe invariant, there exists a unique $\tilde{g}=\varphi^{\frac{4}{n-2}}g$ with zero scalar curvature and $Vol_{\tilde{g}}\left(M\right)=1$. The conditions (\ref{eqNC1}) and (\ref{eqNC2}) of theorem \ref{mainthm} can be rewritten in terms of $g$ as follows~: $\max_M f>0$ and $\int_M f\varphi^{2^\star}\, dv_g<0$.

\section{Proof of the theorem}

We assume in the following that $(M,g)$ has null scalar curvature and that $Vol_g\left(M\right)=1$. We 
let $f\in C^\infty\left(M\right)$, $f\not\equiv 0$, be such that (\ref{eqNC1}) and (\ref{eqNC2}) hold. We first 
recall, and give a quick proof of, the following classical result which deals with the subcritical equation~:

\begin{claim}
\label{propsouscritique}
Let $(M,g)$ be a smooth compact Riemannian manifold of dimension $n\ge 3$ 
with null scalar curvature. Let $f\in C^\infty\left(M\right)$, $f\not\equiv 0$, be such that (\ref{eqNC1}) and (\ref{eqNC2}) hold. For any $2<q<2^\star$, there exists $u_q\in C^\infty\left(M\right)$, $u_q>0$, 
such that 
$$\int_M \left\vert \nabla u_q\right\vert_g^2\, dv_g = \lambda_q\hskip.1cm,\hskip.2cm u_q\in 
{\mathcal H}_q$$
where 
$${\mathcal H}_q = \left\{u\in H_1^2\left(M\right)\hbox{ s.t. }\int_M f\vert u\vert^q\, dv_g=1\right\}$$
and 
$$\lambda_q = \inf_{u\in {\mathcal H}_q} \int_M \left\vert \nabla u\right\vert_g^2\, dv_g\hskip.1cm.$$
Moreover, $u_q$ satisfies that 
$$\Delta_g u_q = \lambda_q f u_q^{q-1}\hskip.1cm.$$
\end{claim}

\medskip {\it Proof} - It is rather standard and we just outline the proof. First, under the assumptions 
made on $f$, it is easily checked that ${\mathcal H}_q$ is not empty. Let $\left(u_i\right)$ be a sequence 
of functions in ${\mathcal H}_q$ such that 
$$\int_M \left\vert \nabla u_i\right\vert_g^2\, dv_g \to \lambda_q$$
as $i\to +\infty$. If we prove that $\left(u_i\right)$ is bounded in $H_1^2\left(M\right)$, then we can extract a subsequence which converges to a solution of our problem. It is straightforward since the embedding of $H_1^2\left(M\right)$ into $L^q\left(M\right)$ is compact for $q<2^\star$. 
Moreover, our minimizer may be choosen nonnegative and the result will follow from the Euler-Lagrange equation satisfied by the minimizer, from standard regularity theory and from the maximum principle. Thus it remains to prove that $\left(u_i\right)$ is bounded in $H_1^2\left(M\right)$. Since $\left(\left\vert \nabla u_i\right\vert_g\right)$ is bounded in $L^2\left(M\right)$, we just have to prove that $\left(u_i\right)$ is bounded in $L^2\left(M\right)$. Assume by 
contradiction that $\Vert u_i\Vert_2\to +\infty$ as $i\to +\infty$. Let $v_i=\frac{u_i}{\Vert u_i\Vert_2}$. Then 
it is clear that $\left(v_i\right)$ is bounded in $H_1^2\left(M\right)$ and that $\left\Vert \nabla v_i\right\Vert_2\to 0$ as $i\to +\infty$. Thus $v_i\to v_\infty$ strongly in $H_1^2\left(M\right)$ as 
$i\to +\infty$ where $v_\infty\equiv 1$. Remember that $Vol_g\left(M\right)=1$. We thus obtain that 
$$0=\lim_{i\to +\infty} \int_M f\left\vert v_i\right\vert^q\, dv_g = \int_M f\left\vert v_\infty\right\vert^q\, dv_g
=\int_M f\, dv_g\hskip.1cm,$$
which is a contradiction since $\int_M f\, dv_g<0$. This ends the proof of the claim. \hfill $\diamondsuit$ 

\medskip Let us now consider the solution $u_q$ of 
\begin{equation}
\label{equq}
\Delta_g u_q = \lambda_q f u_q^{q-1}\hskip.1cm,\hskip.2cm \int_M fu_q^q\, dv_g=1\hskip.1cm,
\end{equation}
given by claim \ref{propsouscritique}. Note that $\lambda_q> 0$ since, otherwise, $u_q$ would be constant and the constraint $\int_M fu_q^q\, dv_g=1$ could not be satisfied since $\int_M f\, dv_g<0$. If $\left(u_q\right)$ is bounded in $L^\infty\left(M\right)$, then standard elliptic theory gives that, after passing to a subsequence, $u_q\to u_0$ in $C^2\left(M\right)$ as $q\to 2^\star$ where 
$u_0$ satisfies that 
$$\Delta_g u_0 = \lambda_{2^\star} f u_0^{2^\star-1}\hskip.1cm,\hskip.2cm \int_M fu_0^{2^\star}\, dv_g=1\hskip.1cm.$$
Indeed, it is easy to check that
\begin{equation}
\label{eqlimitlambdaq}
\lim_{q\to 2^\star} \lambda_q = \lambda_{2^\star}\hskip.1cm.
\end{equation}
Then it is clear that 
$\lambda_{2^\star}> 0$ and that $\lambda_{2^\star}^{\frac{4}{n-2}}u_0$ is a smooth positive solution of (\ref{eqprescribedscalarcurvature}) thanks to the maximum principle and standard regularity theory. 
Thus, in this case, the theorem is proved. We assume now by contradiction that 
\begin{equation}\label{eqassumption}
u_q\left(x_q\right)=\left\Vert u_q\right\Vert_{\infty}\to +\infty\hskip.2cm\hbox{as }q\to 2^\star\hskip.1cm.
\end{equation}
We prove the following~: 

\begin{claim}
\label{proprescaling}
Assume that (\ref{eqassumption}) holds. Then we have that, after passing to a subsequence,  
$$u_q\left(x_q\right)^{-1}u_q\left(\exp_{x_q}\left(\mu_q x\right)\right)\to 
\left(1+\vert x\vert^2\right)^{1-\frac{n}{2}}$$
in $C^2_{loc}\left(\rn\right)$ as $q\to 2^\star$ where 
$$\mu_q = \left(\frac{n(n-2)}{\lambda_q f\left(x_q\right)}\right)^{\frac{1}{2}} u_q\left(x_q\right)^{-\frac{q-2}{2}} \hskip.1cm.$$
\end{claim}

\medskip {\it Proof} - Let us first remark that $\left(u_q\right)$ is bounded in $H_1^2\left(M\right)$. Since 
$\left\Vert \nabla u_q\right\Vert_2^2 \to \lambda_{2^\star}$ as $q\to 2^\star$, we only need to prove that 
$\left(u_q\right)$ is bounded in $L^2\left(M\right)$. Assume on the contrary that $\left\Vert u_q\right\Vert_2\to +\infty$ as $q\to 2^\star$. We then easily obtain that 
$$\frac{u_q}{\left\Vert u_q\right\Vert_2}\to 1$$ 
strongly in $H_1^2\left(M\right)$ as $q\to 2^\star$. This is clearly in contradiction with the fact that 
$\int_M fu_q^q\, dv_g=1$ since $\int_M f\, dv_g<0$. Thus $\left(u_q\right)$ is bounded in $H_1^2\left(M\right)$. As a consequence of Sobolev's embeddings, we know that there exists some $\Lambda>0$ such that 
\begin{equation}\label{eqbounded}
\int_M u_q^q\, dv_g \le \Lambda\hskip.1cm.
\end{equation}
We let now $\delta>0$ small enough and we set for $x\in B_0\left(\delta u_q\left(x_q\right)^{\frac{q-2}{2}}\right)$
\begin{equation}
\label{eqdefrescaling}
\begin{split}
&v_q\left(x\right)=  u_q\left(x_q\right)^{-1}u_q\left(\exp_{x_q}\left(u_q\left(x_q\right)^{-\frac{q-2}{2}}x\right)\right)\hskip.1cm,\\
&g_q(x)=\exp_{x_q}^\star\, g \left(\exp_{x_q}\left(u_q\left(x_q\right)^{-\frac{q-2}{2}}x\right)\right)\hskip.1cm,\\
&f_q\left(x\right)= f\left(\exp_{x_q}\left(u_q\left(x_q\right)^{-\frac{q-2}{2}}x\right)\right)\hskip.1cm.
\end{split}
\end{equation}
It is clear that $g_q\to \xi$, the Euclidean metric, in $C^2_{loc}\left(\rn\right)$ as $q\to 2^\star$, that 
$f_q\to f(x_0)$ in $C^2_{loc}\left(\rn\right)$ as $q\to 2^\star$ where, after passing to a subsequence, 
$x_0=\lim_{q\to 2^\star}x_q$ and that $0\le v_q\le 1$ in $B_0\left(\delta u_q\left(x_q\right)^{\frac{q-2}{2}}\right)$ with $v_q(0)=1$. Moreover, we have that 
$$\Delta_{g_q}v_q = \lambda_q f_q v_q^{q-1}$$
in $B_0\left(\delta u_q\left(x_q\right)^{\frac{q-2}{2}}\right)$. Standard elliptic theory gives then that 
$$v_q\to v_0 \hskip.1cm\hbox{ in }C^2_{loc}\left(\rn\right)\hskip.2cm\hbox{as }q\to 2^\star$$
where $0\le v_0\le 1$, $v_0(0)=1$ and 
$$\Delta_\xi v_0 = \lambda_{2^\star} f(x_0)v_0^{2^\star-1}\hskip.1cm.$$
Since $0$ is a point of maximum of $v_0$, $\Delta_\xi v_0(0)\ge 0$ and thus $f(x_0)\ge 0$. Assume by contradiction that 
\begin{equation}\label{eqfx0equal0}
\lambda_{2^\star} f\left(x_0\right)=0\hskip.1cm.
\end{equation}
Then $v_0$ is a bounded harmonic function in 
$\rn$ with $v_0(0)=1$ and thus $v_0\equiv 1$. 
Now we write thanks to (\ref{eqbounded}) that 
\begincal 
\Lambda &\ge& \int_{B_{x_q}\left(Ru_q\left(x_q\right)^{-\frac{q-2}{2}}\right)} u_q^q\, dv_g \\
&\ge & u_q\left(x_q\right)^{\left(2^\star-q\right)\frac{n-2}{2}}\int_{B_0\left(R\right)} v_q^q\, dv_{g_q} \\
&\ge & \int_{B_0\left(R\right)} v_0^{2^\star}\, dv_\xi+o(1)\hskip.1cm.
\fincal 
This clearly proves that $v_0\in L^{2^\star}\left(\rn\right)$ and thus that $v_0\not\equiv 1$. In particular, (\ref{eqfx0equal0}) can not hold and 
$\lambda_{2^\star}f\left(x_0\right)>0$. It is then clear that 
$$u_q\left(x_q\right)^{-1}u_q\left(\exp_{x_q}\left(\mu_q x\right)\right)\to U_0$$
in $C^2_{loc}\left(\rn\right)$ as $q\to 2^\star$ where 
$\mu_q$ is as in the proposition and $U_0$ satisfies that $0\le U_0\le 1$, $U_0(0)=1$ and 
$$\Delta_\xi U_0= n(n-2) U_0^{2^\star-1}\hskip.1cm.$$
By the classification result of Caffarelli-Gidas-Spruck \cite{CGS}, we know that 
$$U_0(x)=\left(1+\vert x\vert^2\right)^{1-\frac{n}{2}}\hskip.1cm.$$
This ends the proof of the claim. \hfill $\diamondsuit$

\medskip We next prove the following~:

\begin{claim}
\label{propweakconvergence}
Under the assumption (\ref{eqassumption}), we have that, after passing to a subsequence,  $u_q\rightharpoonup 0$ in $H_1^2\left(M\right)$ as $q\to 2^\star$ and that 
$$\lim_{q\to 2^\star} \lambda_q= \lambda_{2^\star}= \left(\max_M f\right)^{-\frac{2}{2^\star}}K(n,2)^{-2}$$
where 
$$K(n,2)^2 = \frac{4}{n(n-2)} \omega_n^{-\frac{2}{n}}$$
with $\omega_n$ the volume of the unit $n$-sphere in ${\mathbb R}^{n+1}$. 
Moreover, we have that 
$$\lim_{R\to +\infty}\lim_{q\to 2^\star} \int_{B_{x_q}\left(R\mu_q\right)}fu_q^q\, dv_g =1\hskip.1cm,$$
that $f(x_0)= \max_M f$ and that 
$$\lim_{q\to 2^\star} u_q\left(x_q\right)^{2^\star-q}=1\hskip.1cm.$$
\end{claim}

\medskip {\it Proof} - Thanks to claim \ref{proprescaling}, we can write that 
\begincal 
&&\int_{B_{x_q}\left(R\mu_q\right)}fu_q^q\, dv_g\\
&& = 
f\left(x_0\right)^{1-\frac{n}{2}}\left(n(n-2)\right)^{\frac{n}{2}} \lambda_q^{-\frac{n}{2}} 
u_q\left(x_q\right)^{\left(2^\star-q\right)\frac{n-2}{2}}\left(\int_{B_0(R)} U_0^{2^\star}\, dx +o(1)\right)\hskip.1cm.
\fincal
As shown during the proof of claim \ref{proprescaling}, we have that 
$\left(u_q\right)$ is bounded in $H_1^2\left(M\right)$. Thus, by Sobolev's embeddings, $\left(u_q\right)$ is uniformly bounded in $L^q\left(M\right)$. We deduce that, after passing to a subsequence, 
\begin{equation}
\label{eqproof2.3.1}
u_q\left(x_q\right)^{\left(2^\star-q\right)\frac{n-2}{2}}\to A
\end{equation}
as $q\to 2^\star$ for some $A\ge 1$. It is also easily checked that 
$$\lim_{q\to 2^\star} \lambda_q = \lambda_{2^\star}$$
and that 
$$\int_{\rn}U_0^{2^\star}\, dx = \bigl(n(n-2)\bigr)^{-\frac{n}{2}}K(n,2)^{-n}\hskip.1cm.$$ 
This leads to 
\begin{equation}
\label{eqproof2.3.2}
\lim_{R\to +\infty}\lim_{q\to 2^\star} \int_{B_{x_q}\left(R\mu_q\right)}fu_q^q\, dv_g
=  f\left(x_0\right)^{1-\frac{n}{2}}  \lambda_{2^\star}^{-\frac{n}{2}}K(n,2)^{-n}A\hskip.1cm.
\end{equation}
Since $\left(u_q\right)$ is bounded in $H_1^2\left(M\right)$, after passing to a subsequence, 
$u_q\rightharpoonup u_0$ weakly in $H_1^2\left(M\right)$ as $q\to 2^\star$. It is clear thanks to 
(\ref{eqproof2.3.2}) that the convergence can not be strong. We then write that 
\begin{equation}
\label{eqproof2.3.3}
\int_M f\left\vert u_q-u_0\right\vert^q\, dv_g = 1-\int_M f\left\vert u_0\right\vert^q\, dv_g +o(1)
\end{equation}
and that 
\begin{equation}\label{eqproof2.3.4}
\int_M \left\vert \nabla\left(u_q-u_0\right)\right\vert_g^2 = \lambda_{2^\star} - \int_M \left\vert \nabla u_0\right\vert_g^2\, dv_g+o(1)\hskip.1cm.
\end{equation}
Equation (\ref{eqproof2.3.4}) is trivial while equation (\ref{eqproof2.3.3}) is easy to obtain (see for instance \cite{BrezisLieb}). 
By the definition of $\lambda_q$, we have that 
$$\lambda_q \left\vert\int_M f\left\vert u_0\right\vert^q\, dv_g\right\vert^{\frac{2}{q}-1}\int_M f\left\vert u_0\right\vert^q\, dv_g \le 
\left\Vert \nabla u_0\right\Vert_2^2$$
and that 
$$\lambda_q \left\vert\int_M f\left\vert u_q-u_0\right\vert^q\, dv_g\right\vert^{\frac{2}{q}-1}\int_M f\left\vert u_q-u_0\right\vert^q\, dv_g \le 
\left\Vert \nabla \left(u_q-u_0\right)\right\Vert_2^2\hskip.1cm.$$
Passing to the limit in these two inequalities, using (\ref{eqproof2.3.3}) and (\ref{eqproof2.3.4}), we get that 
\begin{equation}
\label{eqproof2.3.5}
\lambda_{2^\star}\left\vert X\right\vert^{\frac{2}{2^\star}-1} X\le \left\Vert \nabla u_0\right\Vert_2^2
\end{equation}
and that
\begin{equation}
\label{eqproof2.3.6}
\lambda_{2^\star} \left\vert 1-X\right\vert^{\frac{2}{2^\star}-1} \left(1-X\right)\le \lambda_{2^\star}-\left\Vert \nabla u_0\right\Vert_2^2
\end{equation}
where 
$$X=\int_M f\left\vert u_0\right\vert^{2^\star}\, dv_g\hskip.1cm.$$ 
Since $\left\Vert \nabla u_0\right\Vert_2^2\le \lambda_{2^\star}$ thanks to (\ref{eqproof2.3.4}), the first inequality ensures that $X\le 1$ while the second inequality ensures that $X\ge 0$. Next we sum (\ref{eqproof2.3.5}) and (\ref{eqproof2.3.6}) to get that 
$$\left\vert X\right\vert^{\frac{2}{2^\star}-1} X + \left\vert 1-X\right\vert^{\frac{2}{2^\star}-1} \left(1-X\right)\le 1\hskip.1cm.$$
Since $0\le X\le 1$, this is possible if and only if $X=0$ or $X=1$. But, if $X=1$, (\ref{eqproof2.3.4}) and 
(\ref{eqproof2.3.5}) imply that $u_q\to u_0$ stronly in $H_1^2\left(M\right)$. As already said, this can not happen. Thus $X=0$ and (\ref{eqproof2.3.6}) implies in turn that $\left\Vert \nabla u_0\right\Vert_2=0$ so that $u_0\equiv Cst$. Since $X=\int_M f\left\vert u_0\right\vert^{2^\star}\, dx=0$ and $\int_M f\, dv_g<0$, this clearly implies that $u_0\equiv 0$. This proves the first part of the claim. 
Thanks to the work of Hebey-Vaugon \cite{HebeyVaugon1, HebeyVaugon2}, we know that there exists $B>0$ such that 
$$\left\Vert u_q\right\Vert_{2^\star}^2\le K(n,2)^2 \left\Vert \nabla u_q\right\Vert_2^2 
+ B \left\Vert u_q\right\Vert_{2}^2$$
for all $q<2^\star$ where $K(n,2)$ is as in the statement of the claim. 
Since the embedding of $H_1^2\left(M\right)$ into $L^2\left(M\right)$ is compact, we have that 
$\left\Vert u_q\right\Vert_{2}\to 0$ as $q\to 2^\star$ so that we get that 
$$\left\Vert u_q\right\Vert_{2^\star}^2\le K(n,2)^2 \left\Vert \nabla u_q\right\Vert_2^2 
+o(1)\hskip.1cm.$$
In particular, we deduce that 
\begincal 
1+\int_M f^{-}u_q^q\, dv_g&=& \int_M f^{+}u_q^q\, dv_g\\
&\le & \left(\max_M f\right) \left\Vert u_q\right\Vert_{2^\star}^q Vol_g\left(M\right)^{\frac{\left(2^\star-q\right)}{2^\star}}\\
&\le &  \left(\max_M f\right)K(n,2)^q \left\Vert \nabla u_q\right\Vert_2^q +o(1)\\
&\le &  \left(\max_M f\right)^{\frac{2}{2^\star}}K(n,2)^{2^\star}\lambda_{2^\star}^{\frac{2^\star}{2}} +o(1)
\fincal
where $f^{+}=\frac{1}{2}\left(f+\vert f\vert\right)$ and $f^{-}=\frac{1}{2}\left(\vert f\vert-f\right)$. 
This leads to 
$$\lambda_{2^\star}\ge \left(\max_M f\right)^{-\frac{2}{2^\star}}K(n,2)^{-2}\hskip.1cm.$$
Standard test-functions computations, we refer to \cite{EscobarSchoen} for instance, show that 
the reverse inequality also holds. Thus we have that 
\begin{equation}\label{eqlambda2star}
\lambda_{2^\star}= \left(\max_M f\right)^{-\frac{2}{2^\star}}K(n,2)^{-2}\hskip.1cm.
\end{equation}
Then we can also deduce that 
\begin{equation}\label{eqpartienegative}
\int_M f^{-}u_q^q\, dv_g \to 0\hskip.2cm\hbox{as }q\to 2^\star\hskip.1cm.
\end{equation}
Let us come back to (\ref{eqproof2.3.2}) with (\ref{eqlambda2star}). We obtain that 
$$\lim_{R\to +\infty}\lim_{q\to 2^\star} \int_{B_{x_q}\left(R\mu_q\right)}fu_q^q\, dv_g
= \left(\frac{\max_M f}{f\left(x_0\right)}\right)^{\frac{n-2}{2}} A\hskip.1cm.$$
Since $\int_M fu_q^q\, dv_g=1$ and thanks to (\ref{eqpartienegative}), we deduce that $f(x_0)=\max_M f$ and that $A=1$. This ends the 
proof of the claim. \hfill $\diamondsuit$

\medskip We now prove the following weak pointwise estimates on $u_q$~:

\begin{claim}
\label{propweakestimate}
There exists $C>0$ such that 
$$d_g\left(x_q,x\right)^{\frac{2}{q-2}}u_q(x)\le C$$
for all $x\in M$ and all $q<2^\star$. Moreover, we have that 
$$\lim_{R\to +\infty}\lim_{q\to 2^\star} \sup_{M\setminus B_{x_q}\left(R\mu_q\right)}d_g\left(x_q,x\right)^{\frac{2}{q-2}}u_q(x)=0\hskip.1cm.$$
\end{claim}

\medskip {\it Proof} - The proof is rather standard and follows the lines of \cite{DruetMathAnnalen}. We also refer to \cite{DHRBook} for this kind of argument. We briefly sketch the proof in the following. 
Let us set 
$$w_q(x)=d_g\left(x_q,x\right)^{\frac{2}{q-2}}u_q(x)$$
and let $y_q\in M$ be such that 
$$w_q\left(y_q\right)=\max_M w_q\hskip.1cm.$$
Assume by contradiction that 
\begin{equation}
\label{eqproof2.4.1}
w_q\left(y_q\right)\to +\infty \hskip.2cm\hbox{as }q\to 2^\star\hskip.1cm.
\end{equation}
Since $M$ is compact, it is clear that $u_q\left(y_q\right)\to +\infty$ as $q\to 2^\star$. 
We let now $\delta>0$ small enough and we set for $x\in B_0\left(\delta u_q\left(y_q\right)^{\frac{q-2}{2}}\right)$
\begin{equation}
\label{eqproof2.4.2}
\begin{split}
&v_q\left(x\right)=  u_q\left(y_q\right)^{-1}u_q\left(\exp_{y_q}\left(u_q\left(y_q\right)^{-\frac{q-2}{2}}x\right)\right)\hskip.1cm,\\
&g_q(x)=\exp_{y_q}^\star\, g \left(\exp_{y_q}\left(u_q\left(y_q\right)^{-\frac{q-2}{2}}x\right)\right)\hskip.1cm,\\
&f_q\left(x\right)= f\left(\exp_{y_q}\left(u_q\left(y_q\right)^{-\frac{q-2}{2}}x\right)\right)\hskip.1cm.
\end{split}
\end{equation}
It is clear that $g_q\to \xi$ and that 
$f_q\to f(y_0)$ in $C^2_{loc}\left(\rn\right)$ as $q\to 2^\star$ where, after passing to a subsequence, 
$y_0=\lim_{q\to 2^\star}y_q$. Moreover, we have that 
$$\Delta_{g_q}v_q = \lambda_q f_q v_q^{q-1}$$
in $B_0\left(\delta u_q\left(y_q\right)^{\frac{q-2}{2}}\right)$. Let now $x\in \rn$ and let us 
write thanks to the choice of $y_q$ we made that 
$$d_g \left(x_q, \exp_{y_q}\left(u_q\left(y_q\right)^{-\frac{q-2}{2}}x\right)\right)v_q(x)^{\frac{q-2}{2}} \le 
d_g\left(x_q,y_q\right)\hskip.1cm.$$
This leads to 
$$\left(d_g\left(x_q,y_q\right)- \vert x\vert u_q\left(y_q\right)^{-\frac{q-2}{2}}\right)v_q(x)^{\frac{q-2}{2}} 
\le d_g\left(x_q,y_q\right)$$
which can be written as 
$$\left(1- \vert x\vert w_q\left(y_q\right)^{-\frac{q-2}{2}}\right)v_q(x)^{\frac{q-2}{2}} \le 1\hskip.1cm.$$
We deduce from (\ref{eqproof2.4.1}) that 
$$v_q\left(x\right)\le 1+o(1)\hskip.1cm.$$
Standard elliptic theory gives then that 
$$v_q\to v_0 \hskip.1cm\hbox{ in }C^2_{loc}\left(\rn\right)\hskip.2cm\hbox{as }q\to 2^\star$$
where $0\le v_0\le 1$, $v_0(0)=1$ and 
$$\Delta_\xi v_0 = \lambda_{2^\star} f(y_0)v_0^{2^\star-1}\hskip.1cm.$$
Mimicking the proof of claim \ref{proprescaling}, one proves that $f(y_0)>0$. We then write since 
$\left\Vert u_q\right\Vert_\infty^{2^\star-q}\to 1$ as $q\to 2^\star$, see claim \ref{propweakconvergence}, that 
\begin{equation}
\label{eqproof2.4.3}
\lim_{q\to 2^\star}\int_{B_{y_q}\left(Ru_q\left(y_q\right)^{-\frac{q-2}{2}}\right)} fu_q^q\, dv_g
= f(y_0) \int_{B_0\left(R\right)} v_0^{2^\star}\, dx >0\hskip.1cm.
\end{equation}
Note that (\ref{eqproof2.4.1}) implies that 
$$\frac{d_g\left(x_q,y_q\right)}{u_q\left(y_q\right)^{-\frac{q-2}{2}}}\to +\infty$$
and that 
$$\frac{d_g\left(x_q,y_q\right)}{\mu_q}\to +\infty$$
as $q\to 2^\star$ thanks to the definition of $\mu_q$ and of $x_q$. Thus we can write thanks to (\ref{eqpartienegative}) that for any $R>0$, 
$$\int_M fu_q^q\, dv_g \ge \int_{B_{x_q}\left(R\mu_q\right)} fu_q^q\, dv_g
+\int_{B_{y_q}\left(Ru_q\left(y_q\right)^{-\frac{q-2}{2}}\right)} fu_q^q\, dv_g+o(1)$$
for $q$ close to $2^\star$ which leads to a contradiction thanks to claim \ref{propweakconvergence} and to (\ref{eqproof2.4.3}) since $\int_M fu_q^q\, dv_g=1$. Thus we have proved the first part of the claim. 
Note that, as a consequence of this estimate, we know that 
$\left(u_q\right)$ is uniformly bounded in any compact subset of $M\setminus \left\{x_0\right\}$. Standard  elliptic theory then gives since $u_q\rightharpoonup 0$ in $H_1^2\left(M\right)$ that 
\begin{equation}
\label{eqproof2.4.4}
u_q\to 0 \hskip.2cm\hbox{in }C^2_{loc}\bigl(M\setminus \left\{x_0\right\}\bigr)\hskip.2cm\hbox{as }q\to 2^\star\hskip.1cm.
\end{equation}
Let us now assume by contradiction that the second estimate 
of the proposition does not hold. In other words, let us assume that there exists $z_q\in M$ such that 
\begin{equation}\label{eqproof2.4.5}
\frac{d_g\left(x_q,z_q\right)}{\mu_q}\to +\infty \hskip.2cm\hbox{as }q\to 2^\star
\end{equation}
and that 
\begin{equation}
\label{eqproof2.4.6}
w_q\left(z_q\right)\ge \eps_0
\end{equation}
for some $\eps_0>0$. Note that (\ref{eqproof2.4.6}) together with (\ref{eqproof2.4.4}) implies that 
$d_g\left(x_q,z_q\right)\to 0$ as $q\to 2^\star$ and thus that $u_q\left(z_q\right)\to +\infty$ as $q\to 2^\star$. One can then check that, after passing to a subsequence, 
$$u_q\left(z_q\right)^{-1}u_q\left(\exp_{z_q}\left(u_q\left(z_q\right)^{-\frac{q-2}{2}}x\right)\right)\to 
V_0$$
in $C^2_{loc}\left(\rn\setminus \left\{P\right\}\right)$ where 
$$P=\lim_{q\to 2^\star} u_q\left(z_q\right)^{\frac{q-2}{2}}\exp_{z_q}^{-1}\left(x_q\right)$$
which does exist and satisfies $\eps_0^{\frac{2}{2^\star-2}}\le \vert P\vert \le C^{\frac{2}{2^\star-2}}$ 
where $C$ is the constant involved in the first estimate of the claim. Moreover, $V_0$ satisfies that 
$V_0(0)=1$ and that 
$$\Delta_\xi V_0 = \lambda_{2^\star}f(x_0) V_0^{2^\star-1}\hskip.2cm\hbox{in }\rn\setminus \left\{P\right\}\hskip.1cm.$$
This implies in particular that 
$$\lim_{q\to 2^\star}\int_{B_{z_q}\left(\frac{1}{2}\eps_0^{\frac{2}{2^\star-2}}u_q\left(z_q\right)^{-\frac{q-2}{2}}\right)}fu_q^q\, dv_g=f\left(x_0\right)\int_{B_0\left(\frac{1}{2}\eps_0^{\frac{2}{2^\star-2}}\right)} V_0^{2^\star}\, dx>0\hskip.1cm.$$
Using (\ref{eqproof2.4.5}) and (\ref{eqproof2.4.6}), one proves that, for any $R>0$, the balls $B_{x_q}\left(R\mu_q\right)$ and 
$B_{z_q}\left(\frac{1}{2}\eps_0^{\frac{2}{2^\star-2}}u_q\left(z_q\right)^{-\frac{q-2}{2}}\right)$ are disjoint for $q$ close to $2^\star$. And one obtains a contradiction as in the first part of the proof. This ends the proof of the claim.\hfill $\diamondsuit$

\medskip We transform now this weak pointwise estimate into an almost optimal pointwise estimate. Once again, we refer to \cite{DHRBook} for the general scheme of this kind of proof. 

\begin{claim}
\label{proppointwiseestimate} 
For any $0<\eps<1$, there exists $C_\eps>0$ such that 
$$u_q(x)\le C_\eps\left(\mu_q^{\frac{n-2}{2}\left(1-2\eps\right)}d_g\left(x_q,x\right)^{\left(2-n\right)\left(1-\eps\right)}
+ \eta_q d_g\left(x_q,x\right)^{\left(2-n\right)\eps}\right)$$
for all $x\in M\setminus\left\{x_q\right\}$ and all $q<2^\star$. Here, $\eta_q$ is defined by 
$$\eta_q = \sup_{M\setminus B_{x_q}\left(\delta\right)} u_q$$
where $\delta>0$ is fixed small enough. 
\end{claim}

\medskip {\it Proof} - We let $G$ be the Green function for the Laplacian on $M$ 
normalized such that 
$$\int_M G\left(x,y\right)\, dv_g(y)=0$$
for all $x\in M$. We let $G_q(x)=G\left(x_q,x\right)$ and we let $\delta>0$ be such that 
\begin{equation}\label{eqproof2.5.1}
d_g\left(x_q,x\right)\frac{\left\vert \nabla G_q(x)\right\vert}{G_q(x)} \ge \frac{n-2}{2}
\hskip.2cm\hbox{ and } G_q(x)\ge 1
\end{equation}
for $x\in B_x\left(2\delta\right)$. For the construction of and estimates on the Green function, we refer the reader to the Appendix of \cite{DHRBook}. We fix $0<\eps<1$ and we choose $R_\eps>0$ such that 
\begin{equation}\label{eqproof2.5.2}
\lambda_q d_g\left(x_q,x\right)^2u_q\left(x\right)^{q-2} \le \frac{\eps\left(1-\eps\right)}{f(x_0)} \frac{\left(n-2\right)^2}{8} \hskip.2cm\hbox{in }M\setminus B_{x_q}\left(R_\eps\mu_q\right)
\end{equation}
for $q$ close to $2^\star$ thanks to claim \ref{propweakestimate}. 
We let $L_q$ be the linear operator defined by 
$$L_q \varphi = \Delta_g \varphi -\lambda_q f u_q^{q-2}\varphi\hskip.1cm.$$
This operator satisfies the maximum principle on $M$ in the following sense (see \cite{BNV})~: if $\Omega$ is a smooth subset of $M$ and if $\varphi$ and $\psi$ are two smooth functions in a neighbourhood of $\Omega$ such that 
$$L_q \varphi\ge L_q\psi\hskip.1cm\hbox{ in }\Omega$$
and 
$$\varphi\ge \psi\hskip.1cm\hbox{ on }\partial\Omega\hskip.1cm,$$
then $\varphi\ge \psi$ in $\Omega$. We note that $L_q u_q=0$ in $M$. Simple computations lead thanks 
to (\ref{eqproof2.5.1}) and (\ref{eqproof2.5.2}) to 
\begincal
\frac{L_q G_q^{\nu}}{G_q^\nu} &=& \nu\left(1-\nu\right) \frac{\left\vert \nabla G_q\right\vert^2}{G_q^2}
-\lambda_q fu_q^{q-2}\\
&\ge &  d_g\left(x_q,x\right)^{-2}\left(\nu\left(1-\nu\right)\frac{(n-2)^2}{4} 
- \lambda_q f(x_0) d_g\left(x_q,x\right)^2 u_q(x)^{q-2}\right)\\
&\ge & 0
\fincal
in $B_{x_q}\left(2\delta\right)\setminus B_{x_q}\left(R_\eps\mu_q\right)$ for $q$ close to $2^\star$ and 
for $\nu=\eps$ and $\nu=1-\eps$. Using now claim \ref{proprescaling} and standard estimates 
on the Green function, we obtain the existence of some $C_\eps>0$ such that 
$$u_q\le C_\eps \left(\mu_q^{\frac{n-2}{2}\left(1-2\eps\right)} G_q\left(x\right)^{1-\eps}+ \eta_q G_q\left(x\right)^\eps\right)$$ 
on $\partial \left(B_{x_q}\left(2\delta\right)\setminus B_{x_q}\left(R_\eps\mu_q\right)\right)$. Note that we used the fact that $u_q\left(x_q\right)^{2^\star-q}\to 1$ as $q\to 2^\star$, which was proved in claim  \ref{propweakconvergence}. The maximum principle and standard estimates on the Green function permit to conclude the proof of the claim. Note that, outside $B_{x_q}\left(2\delta\right)$, the estimate of the claim is obviously satisfied, up to change $C_\eps$. \hfill $\diamondsuit$

\medskip Remark that, by standard elliptic theory, see for instance \cite{GT}, chapter 8, 
we know that 
$$\frac{u_q}{\eta_q}\to H$$
in $C^2_{loc}\left(M\setminus\left\{x_0\right\}\right)$ as $q\to 2^\star$ where $H$ is a nonzero positive function satisfying that $\Delta_g H=0$ in $M\setminus\left\{x_0\right\}$. Multiplying the equation (\ref{equq}) by $u_q^{1-q}$ and integrating over $M$, we obtain that 
\begincal
\lambda_q\int_M f\, dv_g &=& -\left(q-1\right)\int_M u_q^{-q}\left\vert \nabla u_q\right\vert_g^2\, dv_g\\
&\le & -\left(q-1\right)\eta_q^{2-q}\left(\int_{M\setminus B_{x_0}\left(\delta\right)} H^{-q}\left\vert \nabla H\right\vert^2\, dv_g+o(1)\right)
\fincal
for all $\delta>0$. Since $\eta_q\to 0$ as $q\to 2^\star$ thanks to (\ref{eqproof2.4.4}), this clearly implies that $\nabla H\equiv 0$. Thus, by the definition of $\eta_q$, we get 
that $H\equiv 1$ and we have proved that 
\begin{equation}
\label{eqconvergenceexterieure}
\frac{u_q}{\eta_q}\to 1\hskip.1cm\hbox{ in }C^2_{loc}\bigl(M\setminus\left\{x_0\right\}\bigr)\hskip.1cm\hbox{ as }q\to 2^\star\hskip.1cm.
\end{equation}
We now describe precisely the asymptotic behaviour of $u_q$~:

\begin{claim}
\label{propasymptoticestimates}
For any sequence $\left(y_q\right)$ of points in $M$, we have that 
\begincal
u_q\left(y_q\right) &=& \overline{u_q} + 2^n\frac{\omega_{n-1}}{n}u_q\left(x_q\right)^{-1}\varphi\left(y_q\right)+ o\left(u_q\left(x_q\right)^{-1}\right)\\
&&+ u_q\left(x_q\right) \left(1+\frac{d_g\left(x_q,y_q\right)^2}{\mu_q^2}\right)^{-\frac{n-2}{2}}\bigl(\Phi\left(x_q,y_q\right)+o(1)\bigr)\hskip.1cm.
\fincal
In this equation,  
$$\overline{u_q}= \int_M u_q\, dv_g$$
and $\varphi$ is a solution of 
$$\Delta_g \varphi = \lambda_{2^\star}\left(1-\frac{f}{\bar{f}}\right)$$
with $\bar{\varphi}=0$ and 
$$\Phi(x,y)= \left(n-2\right)\omega_{n-1}d_g\left(x,y\right)^{n-2}G\left(x,y\right)$$
where $G$ is the Green function of the Laplacian on $M$ normalized 
such that 
$$\int_M G\left(x,y\right)\, dv_g(y)=0$$
for all $x\in M$.
\end{claim}

\medskip {\it Proof} - We fix $0<\eps<\frac{2}{n+2}$. We 
first write  by integrating equation (\ref{equq}) over $M$ that 
$$\int_M fu_q^{q-1}\, dv_g=0$$
and thus that 
\begin{equation}\label{eqproof2.6.02}
\int_{\left\{u_q(x)> \alpha \eta_q\right\}} fu_q^{q-1}\, dv_g=
-\int_{\left\{u_q(x)\le \alpha \eta_q\right\}} fu_q^{q-1}\, dv_g
\end{equation}
for all $\alpha>1$. 
Thanks to Lebesgue's dominated convergence theorem and to (\ref{eqconvergenceexterieure}), 
we have that 
\begin{equation}\label{eqproof2.6.03}
\int_{\left\{u_q(x)\le \alpha \eta_q\right\}} fu_q^{q-1}\, dv_g = \eta_q^{q-1}\int_M f\, dv_g +
o\left(\eta_q^{q-1}\right)\hskip.1cm.
\end{equation}
We write also thanks to claim \ref{proprescaling} that, for any $R>0$, 
\begin{equation}\label{eqproof2.6.04}
\begin{split}
\int_{\left\{u_q(x)> \alpha \eta_q\right\}} fu_q^{q-1}\, dv_g=\, &
\int_{B_{x_q}\left(R\mu_q\right)} fu_q^{q-1}\, dv_g \\
&+
\int_{\left\{u_q(x)> \alpha \eta_q\right\}\setminus B_{x_q}\left(R\mu_q\right)}fu_q^{q-1}\, dv_g\hskip.1cm.
\end{split}
\end{equation}
Thanks to claim \ref{proppointwiseestimate}, we have that 
\begincal 
&&\left\vert \int_{\left\{u_q(x)> \alpha \eta_q\right\}\setminus B_{x_q}\left(R\mu_q\right)}fu_q^{q-1}\, dv_g
\right\vert \\
&&\quad \le 
C\mu_q^{\frac{n-2}{2}\left(1-2\eps\right)\left(q-1\right)} \int_{M\setminus B_{x_q}\left(R\mu_q\right)} 
d_g\left(x_q,x\right)^{\left(2-n\right)\left(1-\eps\right)\left(q-1\right)}\, dv_g\\
&&\qquad + C\eta_q^{q-1} \int_{\left\{u_q(x)> \alpha \eta_q\right\}\setminus B_{x_q}\left(R\mu_q\right)}
d_g\left(x_q,x\right)^{\left(2-n\right)\eps\left(q-1\right)}\, dv_g\\
&&\quad \le C R^{n-\left(n+2\right)\left(1-\eps\right)}\mu_q^{\frac{n-2}{2}}+ C \alpha^{-\frac{n-\left(n+2\right)\eps}{\left(n-2\right)\eps}}\eta_q^{q-1}+C \alpha^{\frac{1-2\eps}{\eps}\left(q-1\right)}\mu_q^{\frac{n+2}{2}\left(1-2\eps\right)}
\fincal
where $C$ changes from line to line but is always independent of $R$, $\alpha$ and $q$. 
Here we used the fact that $\eps<\frac{2}{n+2}$ and the fact that $u_q\left(x_q\right)^{2^\star-q}\to 1$ as $q\to 2^\star$, which was proved in claim \ref{propweakconvergence}. We also used the 
following consequence of claim \ref{proppointwiseestimate}~:
\begin{equation}\label{eqproof2.6.01}\begin{split}
&\hbox{Either}\hskip.2cm u_q(x) >\alpha\eta_q \Rightarrow d_g\left(x_q,x\right)^{\left(n-2\right)\eps}\le \frac{2C_\eps}{\alpha}\hskip.1cm,\\
&\hbox{or}\hskip.2cm\eta_q\le \left(\frac{\alpha}{2C_\eps}\right)^{\frac{1-2\eps}{\eps}} \mu_q^{\frac{n-2}{2}\left(1-2\eps\right)}\hskip.1cm.
\end{split} 
\end{equation}
Since $0<\eps<\frac{2}{n+2}$, we deduce from the above estimate 
and the definition of $\mu_q$ that 
\begin{equation}\label{eqproof2.6.05}
\lim_{\alpha\to +\infty}\lim_{R\to +\infty}\lim_{q\to 2^\star} \frac{\int_{\left\{u_q(x)> \alpha \eta_q\right\}\setminus B_{x_q}\left(R\mu_q\right)}fu_q^{q-1}\, dv_g}{u_q\left(x_q\right)^{-1}+\eta_q^{q-1}}=0\hskip.1cm.
\end{equation}
Now, thanks to claims \ref{proprescaling} and \ref{propweakconvergence}, it is easily checked that 
\begin{equation}
\label{eqproof2.6.06}
\lim_{R\to +\infty}\lim_{q\to 2^\star} u_q\left(x_q\right)\int_{B_{x_q}\left(R\mu_q\right)} fu_q^{q-1}\, dv_g 
= \frac{2^n\omega_{n-1}}{n}\hskip.1cm.
\end{equation}
Coming back to (\ref{eqproof2.6.04}) with (\ref{eqproof2.6.05}) and (\ref{eqproof2.6.06}) , we thus get that 
\begin{equation}
\label{eqproof2.6.07}
\int_{\left\{u_q(x)> \alpha \eta_q\right\}} fu_q^{q-1}\, dv_g= \frac{2^n\omega_{n-1}}{n}u_q\left(x_q\right)^{-1}+o\left(u_q\left(x_q\right)^{-1}\right)+o\left(\eta_q^{q-1}\right)\hskip.1cm.
\end{equation}
Coming back to (\ref{eqproof2.6.02}) with (\ref{eqproof2.6.03}) and (\ref{eqproof2.6.07}), we obtain that 
\begin{equation}\label{eqproof2.6.08}
\eta_q^{q-1} = -\frac{2^n\omega_{n-1}}{n\bar{f}}u_q\left(x_q\right)^{-1} + o\left(u_q\left(x_q\right)^{-1}\right)\hskip.1cm.
\end{equation}
We write now thanks to the Green representation formula that 
\begin{equation}\label{eqproof2.6.1}
u_q\left(y_q\right)- \overline{u_q}= \lambda_q\int_M G\left(y_q,x\right) f(x)u_q(x)^{q-1}\, dv_g (x)
\end{equation}
for all sequences $y_q\in M$. We write thanks to Lebesgue's dominated convergence theorem, to (\ref{eqconvergenceexterieure}) and to (\ref{eqproof2.6.08}) that 
\begincal
&&\int_M G\left(y_q,x\right) f(x)u_q(x)^{q-1}{\bf 1}_{\left\{u_q\le \alpha \eta_q\right\}}\, dv_g (x)\\
&&\quad = \eta_q^{q-1}\int_M G\left(y_q,x\right) f(x)\, dv_g + o\left(\eta_q^{q-1}\right)\\
&&\quad =  -\frac{2^n\omega_{n-1}}{n\bar{f}}u_q\left(x_q\right)^{-1}\int_M G\left(y_q,x\right) f(x)\, dv_g
+o\left(u_q\left(x_q\right)^{-1}\right)
\fincal
for all $\alpha>0$. We can write now that 
\begincal 
\int_M G\left(y_q,x\right) f(x)\, dv_g &=& \int_M G\left(y_q,x\right) \left(f(x)-\bar{f}\right)\, dv_g\\
&=& -\bar{f}\frac{1}{\lambda_2^\star} \varphi\left(y_q\right)
\fincal 
so that we obtain that 
\begin{equation}\label{eqproof2.6.2}
\begin{split}
&\lambda_q\int_M G\left(y_q,x\right) f(x)u_q(x)^{q-1}{\bf 1}_{\left\{u_q\le \alpha \eta_q\right\}}\, dv_g (x)\\
&\quad= \frac{2^n\omega_{n-1}}{n}u_q\left(x_q\right)^{-1}\varphi\left(y_q\right)
+o\left(u_q\left(x_q\right)^{-1}\right)\hskip.1cm.
\end{split}
\end{equation}
We assume in the following that $\frac{d_g\left(x_q,y_q\right)}{\mu_q}\to +\infty$ since, otherwise, the estimate of the claim is a straightforward consequence of claim \ref{proprescaling}. We write thanks to claim \ref{proprescaling} that 
\begincal 
&&\int_M G\left(y_q,x\right) f(x)u_q(x)^{q-1}{\bf 1}_{\left\{u_q> \alpha \eta_q\right\}}\, dv_g (x)\\
&&\quad = \int_{B_{x_q}\left(R\mu_q\right)}G\left(y_q,x\right) f(x)u_q(x)^{q-1}\, dv_g (x)\\
&&\qquad + \int_{M\setminus B_{x_q}\left(R\mu_q\right)}G\left(y_q,x\right) f(x)u_q(x)^{q-1}
{\bf 1}_{\left\{u_q> \alpha \eta_q\right\}}\, dv_g (x)\hskip.1cm.
\fincal
Then, thanks to the estimate of claim \ref{proppointwiseestimate}, to (\ref{eqproof2.6.01}) and (\ref{eqproof2.6.08}) - note that (\ref{eqproof2.6.08}) prevents the second situation from happening in (\ref{eqproof2.6.01}), and 
to standard estimates on the Green function, we can write that 
\begincal 
&&\int_{M\setminus B_{x_q}\left(R\mu_q\right)}G\left(y_q,x\right) f(x)u_q(x)^{q-1}
{\bf 1}_{\left\{u_q> \alpha \eta_q\right\}}\, dv_g (x)\\
&&\quad \le C\mu_q^{\frac{n+2}{2}\left(1-2\eps\right)}\int_{M\setminus B_{x_q}\left(R\mu_q\right)} d_g\left(y_q,x\right)^{2-n} d_g\left(x_q, x\right)^{\left(2-n\right)\left(1-\eps\right)\left(q-1\right)}\, dv_g(x)\\
&&\qquad + Cu_q\left(x_q\right)^{-1} \int_{M\setminus B_{x_q}\left(R\mu_q\right)}d_g\left(y_q,x\right)^{2-n} d_g\left(x_q,x\right)^{\left(2-n\right)\eps}{\bf 1}_{\left\{u_q> \alpha \eta_q\right\}}\, dv_g (x)\\
&&\quad \le C R^{n-\left(n+2\right)\left(1-\eps\right)}\mu_q^{\frac{n-2}{2}}d_g\left(x_q,y_q\right)^{2-n}
+C u_q\left(x_q\right)^{-1} \alpha^{-\frac{2-\left(n+2\right)\eps}{\left(n-2\right)\eps}}
\fincal
for some $C>0$ independent of $R$, $\alpha$ and $q$ thanks to the fact that $\eps< \frac{2}{n+2}$. In 
particular, we obtain that 
\begin{equation}
\label{eqproof2.6.3}
\lim_{\alpha\to +\infty}\lim_{R\to +\infty}\lim_{q\to 2^\star}\frac{\int_{M\setminus B_{x_q}\left(R\mu_q\right)}G\left(y_q,x\right) f(x)u_q(x)^{q-1}
{\bf 1}_{\left\{u_q> \alpha \eta_q\right\}}\, dv_g (x)}{u_q\left(x_q\right)^{-1}d_g\left(x_q,y_q\right)^{2-n} }=0\hskip.1cm.
\end{equation}
It is now easily checked thanks to claim \ref{proprescaling} that 
\begin{equation}
\label{eqproof2.6.4}
\begin{split}
&\lim_{R\to +\infty}\lim_{q\to 2^\star} \frac{\int_{B_{x_q}\left(R\mu_q\right)}G\left(y_q,x\right) f(x)u_q(x)^{q-1}\, dv_g (x)}{u_q\left(x_q\right)^{-1}G\left(y_q,x_q\right)}\\
&\quad = \frac{\omega_{n-1}}{n} \left(\frac{n(n-2)}{\lambda_{2^\star}}\right)^{\frac{n}{2}}f\left(x_0\right)^{1-\frac{n}{2}}\hskip.1cm.
\end{split} 
\end{equation}
Coming back to (\ref{eqproof2.6.1}) with (\ref{eqproof2.6.2}), (\ref{eqproof2.6.3}) and (\ref{eqproof2.6.4}), we obtain using claim  
\ref{propweakconvergence} that 
\begincal 
u_q\left(y_q\right) -\overline{u_q}&= &\left(\left(n-2\right)\omega_{n-1}\left(\frac{2^n}{\omega_n f\left(x_0\right)}\right)^{\frac{n-2}{2}}+o(1)\right) 
u_q\left(x_q\right)^{-1}G\left(y_q,x_q\right)\\
&& + \frac{2^n\omega_{n-1}}{n}u_q\left(x_q\right)^{-1}\varphi\left(y_q\right)+ o\left(u_q\left(x_q\right)^{-1}\right)
\fincal
when $\frac{d_g\left(x_q,y_q\right)}{\mu_q}\to +\infty$. Using claim \ref{proprescaling}, 
one ends the proof.\hfill $\diamondsuit$

\medskip We conclude the proof of the theorem. We set 
$$R_q = u_q\left(x_q\right)\left(u_q-\overline{u_q}\right)$$
and we write that 
$$\Delta_g R_q = \lambda_q u_q\left(x_q\right) f u_q^{q-1}\hskip.1cm.$$
Let $K$ be a compact subset of $M\setminus \left\{x_0\right\}$. Thanks to proposition \ref{propasymptoticestimates}, we know that 
$$R_q \to 2^n\frac{\omega_{n-1}}{n}\varphi +\left(n-2\right)\omega_{n-1} 2^{n-2}
f\left(x_0\right)^{-\frac{n-2}{n}}G\left(x_0,\, .\,\right) $$
in $C^0\left(K\right)$ as $q\to 2^\star$. Using also (\ref{eqconvergenceexterieure}) and  (\ref{eqproof2.6.08}), we get that 
$$\Delta_g R_q \to -\frac{2^n \omega_{n-1}}{n}\lambda_{2^\star} \frac{f}{\bar{f}}$$
in $C^1\left(K\right)$ as $q\to 2^\star$. Thus, by standard elliptic theory, the above convergence of $R_q$ holds in $C^2\left(K\right)$ and we can pass to the limit in the equation satisfied by $R_q$ to obtain that 
\begincal
-\frac{2^n \omega_{n-1}}{n}\lambda_{2^\star} \frac{f}{\bar{f}}&=& \frac{2^n\omega_{n-1}}{n}\Delta_g \varphi +\left(n-2\right)\omega_{n-1} 2^{n-2}
f\left(x_0\right)^{-\frac{n-2}{n}}\Delta_g G\left(x_0,\, .\,\right)\\
&=&  \frac{2^n\omega_{n-1}}{n} \lambda_{2^\star}\left(1-\frac{f}{\bar{f}}\right)
\fincal
which is clearly a contradiction. This proves the theorem. 

\medskip {\bf Remarks} : the solution we obtain is a minimizing solution in the sense we have that 
$$\int_M \left\vert \nabla u\right\vert_g^2 \, dv_g= \lambda_{2^\star}(f)$$
and $u\in {\mathcal H}_{2^\star}$ where $\lambda_{2^\star}(f)$ and ${\mathcal H}_{2^\star}(f)$ are as 
in claim \ref{propsouscritique}. We add the dependence in $f$ with respect to the notations of claim \ref{propsouscritique}. Note also that as a consequence of this result, we have that 
\begin{equation}\label{eqremark1}
\lambda_{2^\star}(f)< K(n,2)^{-2} \left(\max_M f\right)^{-\frac{2}{2^\star}}
\end{equation}
for all smooth functions $f$ satisfying (\ref{eqNC1}) and (\ref{eqNC2}), a result which is not obvious 
a priori. Indeed, imagine that, on the contrary, there exists some $f\in C^\infty\left(M\right)$ satisfying 
 (\ref{eqNC1}) and (\ref{eqNC2}) and such that 
\begin{equation}\label{eqremark2}
\lambda_{2^\star}\left(f\right)=K(n,2)^{-2} \left(\max_M f\right)^{-\frac{2}{2^\star}}\hskip.1cm.
\end{equation}
Let us look at the function $\tilde{f} = -\max_M f +2f$ which satisfies (\ref{eqNC1}) and (\ref{eqNC2}) 
and which satisfies moreover that $\max_M \tilde{f}=\max_M f$. 
By test functions computations we know that 
$$\lambda_{2^\star}\left(\tilde{f}\right)\le K(n,2)^{-2} \left(\max_M f\right)^{-\frac{2}{2^\star}}\hskip.1cm.$$
Moreover, thanks to what we proved above, there exists $u\in C^\infty\left(M\right)$, $u>0$, which satisfies 
that 
$$\int_M \tilde{f}u^{2^\star}\, dv_g=1$$ 
and that 
$$\int_M \left\vert \nabla u\right\vert_g^2\, dv_g = \lambda_{2^\star}\left(\tilde{f}\right)\le 
K(n,2)^{-2} \left(\max_M f\right)^{-\frac{2}{2^\star}}\hskip.1cm.$$
Since $f\ge \tilde{f}$ and $f\not\equiv \tilde{f}$, we have that 
$$\int_M fu^{2^\star}\, dv_g >1$$ 
and we can write thanks to the definition of $\lambda_{2^\star}(f)$ that 
$$\lambda_{2^\star}(f)\left(\int_M f u^{2^\star}\, dv_g\right)^{\frac{2}{2^\star}} 
\le \int_M \left\vert \nabla u\right\vert_g^2\, dv_g\le K(n,2)^{-2} \left(\max_M f\right)^{-\frac{2}{2^\star}}
\hskip.1cm.$$
This is clearly in contradiction with (\ref{eqremark2}). Thus the above claim (\ref{eqremark1}) is proved. 

\bigskip {\bf Acknowledgements}~: I would like to thank Emmanuel Hebey and Fr\'ed\'eric Robert 
for their interest in this work and their fruitful comments about this paper.

\bibliographystyle{amsplain}
\bibliography{Refprescribednull}

\end{document}